\documentclass{amsart}

\usepackage{amsmath,amsthm,amsfonts,amssymb,bm,graphicx, mathrsfs}

\usepackage
{hyperref}
\hypersetup{colorlinks=true,citecolor=blue,linkcolor=blue,urlcolor=blue,
pdfstartview=FitH }

\theoremstyle{plain}
\newtheorem{theorem}{Theorem}

\newtheorem{lemma}{Lemma}

\newtheorem{cor}[lemma]{Corollary}

\numberwithin{equation}{section}

\theoremstyle{definition}

\renewcommand{\geq}{\geqslant}
\renewcommand{\leq}{\leqslant}

\newcommand{\changed}[1]{{\color{black} #1}}

\usepackage{calc}
\newsavebox\CBox
\newcommand\hcancel[2][0.5pt]{%
  \changed{\ifmmode\sbox\CBox{$#2$}\else\sbox\CBox{#2}\fi%
  \makebox[0pt][l]{\usebox\CBox}%
  \rule[0.5\ht\CBox-#1/2]{\wd\CBox}{#1}}}

\makeatletter
\DeclareRobustCommand\widecheck[1]{{\mathpalette\@widecheck{#1}}}
\def\@widecheck#1#2{%
    \setbox\z@\hbox{\m@th$#1#2$}%
    \setbox\tw@\hbox{\m@th$#1%
       \widehat{%
          \vrule\@width\z@\@height\ht\z@
          \vrule\@height\z@\@width\wd\z@}$}%
    \dp\tw@-\ht\z@
    \@tempdima\ht\z@ \advance\@tempdima2\ht\tw@ \divide\@tempdima\thr@@
    \setbox\tw@\hbox{%
       \raise\@tempdima\hbox{\scalebox{1}[-1]{\lower\@tempdima\box
\tw@}}}%
    {\ooalign{\box\tw@ \cr \box\z@}}}
\makeatother

\begin{document}

\author{Valentin Blomer}
 
\address{Mathematisches Institut, Bunsenstr. 3-5, 37073 G\"ottingen, Germany} \email{vblomer@math.uni-goettingen.de}
 
 \title{Higher order divisor problems }

\thanks{Author   supported in part  by the  Volkswagen Foundation and an NSF grant while enjoying the hospitality of the Institute for Advanced Study}


\begin{abstract} An asymptotic formula is proved for the $k$-fold divisor function averaged over   homogeneous polynomials  of degree $k$ in $k-1$ variables coming from incomplete norm forms. 
 \end{abstract}


\setcounter{tocdepth}{2}  \maketitle 

\maketitle

\section{Introduction}

The divisor function $\tau  = \tau_2$ and its higher order relatives $\tau_k$ defined by 
\begin{equation}\label{tauk}
\tau_k(n) := |\{(n_1, \ldots, n_k) \in \Bbb{N}^k  \mid n_1 \cdot \ldots \cdot n_k = n\}|, 
\end{equation}
i.e.\ the coefficients of the Dirichlet series for $\zeta^k$,  belong to the most prominent arithmetic functions.  The values $\tau_k(n)$ fluctuate quite considerably as $n$ varies, but the average behaviour is reasonably stable. Their statistical behaviour can be measured in various ways, most interestingly perhaps by considering mean values over sparse sequences. In this paper we consider the values of $\tau_k$ on a thin sequence $\mathcal{S} \subseteq \Bbb{N}$ of logarithmic density $1 - 1/k$, i.e.\   $\mathcal{S}(X) = \mathcal{S} \cap [1, X]$ satisfies
\begin{equation}\label{12}
|\mathcal{S}(X)| \sim \text{const} \cdot X^{1 - 1/k}
\end{equation}
as $X \rightarrow \infty$. One might expect an asymptotic formula of the type
$$\sum_{ n \in \mathcal{S}(X^k)} \tau_k(n) = (C+ o(1))  X^{k-1}  (\log X)^{k-1}$$
for a certain constant $C$ depending on $\mathcal{S}$. The most natural approach to such a result consists in    opening the divisor function $\tau_k$  and  writing
 \begin{equation}\label{13}
\sum_{ n \in \mathcal{S}(X^k)} \tau_k(n) = \sum_{n_2 \cdot \ldots \cdot n_k \leq X^k} |\{ n \in \mathcal{S}(X^k) \mid   n \equiv 0 \, (\text{mod } n_2 \cdot \ldots \cdot n_k) \}|.
\end{equation}
By a symmetry argument (a variation of Dirichlet's hyperbola method) we may essentially assume that
$n_1 \geq n_2 \geq \ldots \geq n_k$ in \eqref{tauk}, so that $n_2 \cdot \ldots \cdot n_k \leq X^{k-1}$.  Thus we need to understand the elements of $\mathcal{S}(X^k)$ in residue classes to moduli of size up to $X^{k-1}$, i.e.\ we need to show level of distribution $1-1/k$. Comparing with \eqref{12}, this is on the edge of what one can expect to be able to prove:  we need to obtain asymptotic information on sets  appearing on the right-hand side of \eqref{13} which generically have   \emph{only a bounded number of elements}.  It can therefore be expected that this leads to a  fairly delicate counting problem. 

We now describe the type of sequences $\mathcal{S}$ we have in mind. Let $K/\Bbb{Q}$ be a  Galois  number field  of  degree $k \geq 3$. Let $\{1=\omega_1, \ldots, \omega_k\}$ be an integral basis of the ring of integers $\mathcal{O}_K$. The associated norm form is given by
$$N_{K/\Bbb{Q}}(x_1 \omega_1  + \ldots x_k \omega_k)  \in \Bbb{Z}[x_1, \ldots x_k].$$
Since $K$ is fixed throughout the paper, we will drop it from the notation and simply write $N$ for the norm. Let 
$$f(\textbf{x}) := N(x_1 \omega_1  + \ldots +  x_{k-1} \omega_{k-1}) \in \Bbb{Z}[x_1, \ldots x_{k-1}]$$
be the incomplete norm form with vanishing last coordinate. For instance, for the biquadratic  field $K = \Bbb{Q}(\sqrt{2} , i)$ with basis $\{1, \sqrt{2}, i, \frac{1}{2}(\sqrt{2}+ \sqrt{2} i)\}$ (see \cite{W}) we have
$$f(x_1, x_2, x_3) = (x_1^2 - 2x_2^2)^2 + x_3^2(2x_1^2 + 4x_2^2 + x_3^2).$$
We denote generally by $\mathcal{R} \subseteq  \Bbb{R}^{k-1}$  a region with piecewise smooth boundary not containing $0 \in \Bbb{R}^{k-1}$ such that any line parallel to the coordinate axes intersects $\mathcal{R}$ in $O(1)$ intervals. For $X > 1$ we write $\mathcal{R}_X := \{X \cdot \textbf{x} \in \Bbb{R}^{k-1} \mid \textbf{x} \in \mathcal{R}\}$. We suppose that $f(\mathcal{R}) \subseteq [1, 2]$ and  are interested in an asymptotic formula for 
$$\mathcal{M}(\mathcal{R}_X) := \sum_{\textbf{x} \in \mathcal{R}_X \cap \Bbb{Z}^{k-1}} \tau_k(f(\textbf{x})).$$
The asymptotic formula will naturally feature the local densities
\begin{equation}\label{varrho}
\varrho(n) := \frac{1}{n^{k-2}} \left|\left\{\textbf{x} \in (\Bbb{Z}/n\Bbb{Z})^{k-1} : n \mid f(\textbf{x})\right\}\right|.
\end{equation} 
Let 
\begin{equation}\label{C}
C := \prod_p \left(\sum_{\nu=0}^{\infty} \frac{\varrho(p^{\nu}) \tau_{k-1}(p^{\nu})}{p^{\nu}}\right)\left(1 - \frac{1}{p}\right)^{k-1}.
\end{equation}
We will see later that this Euler product is absolutely convergent. We 
 have the following main result:
\begin{theorem}\label{thm1} For $\mathcal{R}, K$ as above and $\varepsilon > 0$  we have
$$\mathcal{M}(\mathcal{R}_X) =    \frac{C\cdot    {\rm vol}(\mathcal{R})}{(k-1)!} X^{k-1}   (\log X^k)^{k-1} + O_{K, \mathcal{R}, \varepsilon}\left(X^{k-1} (\log X)^{k-1 - \frac{1}{k-1} + \varepsilon} \right) $$
as $X \rightarrow \infty$. 
\end{theorem}
The case $k=2$ is not particularly interesting. In this case  we have $f(x) = x^2$, so that we are in a situation that can be handled by elementary multiplicative number theory (see e.g.\ \cite[p.\ 89]{D2}). The case $k=3$ was dealt with by Daniel  \cite[p.\ 90]{D2} in somewhat greater generality and with a somewhat weaker error term, using similar techniques as in \cite{D1}. In fact, he proves an asymptotic formula for 
$$\sum_{\substack{\textbf{x} \in \Bbb{Z}^2\\ 0 < |f(\textbf{x})| \leq N}}\tau_3(|f(\textbf{x})|)$$
for an arbitrary irreducible binary cubic form. By  \cite[Lemma 2.1]{HBM} every such form is a norm form, so that his result is in effect a version of Theorem \ref{thm1} for $k=3$ without the assumption that $K$ is Galois.  

A  related problem for $k=4$ of comparable difficulty, although with a polynomial  not coming from an incomplete norm form, was treated in \cite{Ti}, based on methods in \cite{FI}. For larger values of $k$, the author is not aware of a result in this direction. 

We briefly indicate some ingredients of the proof. For simplicity of exposition let us assume for the moment that $K$ has class number 1, i.e.\ unique factorization. The 
 divisibility condition $n \mid  f(\textbf{x})$ can be translated into a (set of) divisibility relations $\textbf{n} \mid \textbf{x}$  in $\mathcal{O}_K$, i.e.\ $\textbf{n} \cdot \textbf{m} = \textbf{x}$. Here we identify integral vectors and integral elements in $K$ using the basis $\{\omega_1, \ldots, \omega_k\}$. As $\textbf{x}$ has vanishing last coordinate, this means that $\textbf{m}$ lies in a $(k-1)$-dimensional  lattice $\Lambda = \Lambda(\textbf{n})$. The key step is to show that in a suitable sense the smallest non-trivial vector in $\Lambda(\textbf{n})$ is   roughly as big as $\text{det}(\Lambda(\textbf{n}))^{1/(k-1)}$ on average over $\textbf{n}$, in other words, $\Lambda(\textbf{n})$ is typically ``well-balanced''.  This analysis  shares some similarities with  recent work of Maynard \cite{Ma} in connection with primes represented by incomplete norm forms. Considerable technical difficulties arise from prime ideals in $\mathcal{O}_K$ of degree $> 1$, which require a   careful setup in Section \ref{sec6}.  The assumption that $K$ is Galois  can probably be dropped, we use it here to simplify the splitting types of rational primes.

It should be emphasized that the main achievement of Theorem \ref{thm1} is to establish an asymptotic formula; upper and lower bounds
$$\mathcal{M}(\mathcal{R}_X)  \asymp  X^{k-1}   (\log X)^{k-1} $$
($k \geq 3$) can be derived from a result of Wolke \cite[Satz 1 \& 2]{Wo} (see  Corollary \ref{wolke1} for the upper bound),  which in turn is an application of Selberg's sieve. In related situations, upper bounds can also be obtained from the work of Nair-Tenenbaum \cite{NT} and Henriot \cite{H}. \\ 

\textbf{Notation:}  The Vinogradov symbols $\ll$ and $\asymp$ have their usual meanings. All implicit constants may depend on the field $K$ (in particular   on $k$) and   on $\mathcal{R}$, and we do not display this dependence.  

The arithmetic function $\omega(n)$ denotes the number of distinct prime divisors   of $n$. For a prime $p$ and natural numbers $\alpha, n$ we write $p^{\alpha} \parallel n$ if $p^{\alpha} \mid n$, but $p^{\alpha+1} \nmid n$. 

 We will often identify vectors $\textbf{x}  = (x_1, \ldots, x_k)\in \Bbb{Z}^k$ and algebraic numbers $\textbf{x}  = x_1 \omega_1 + \ldots + x_k \omega_k \in  \mathcal{O}_K$, and we   denote the corresponding  principal ideal by $(\textbf{x})$.  Multiplication of vectors is defined as the multiplication in $\mathcal{O}_K$; explicitly, if $\omega_i \omega_j = \sum_{r} \alpha_{i, j, r} \omega_r$ for some $\alpha_{i, j, r} \in \Bbb{Z}$, then
$$(x_1, \ldots, x_k) \cdot (y_1, \ldots, y_k) = \Bigl(\sum_{i, j= 1}^n x_i y_j \alpha_{i, j, r}\Bigr)_{1 \leq r \leq k}.$$
The norm of ideals in $\mathcal{O}_K$ is denoted by $N$, whereas $\|. \|$ denotes the Euclidean norm on $\Bbb{R}^{k}$. We embed $\Bbb{Z}^{k-1}$ and $\Bbb{R}^{k-1}$ into $\Bbb{Z}^k$ and $\Bbb{R}^k$ as vectors with vanishing last coordinate. 
We denote by $\textbf{x} \mapsto \textbf{x}^{\sigma}$ an embedding  of $K$ into $\Bbb{R}$ or $\Bbb{C}$, both of which are equipped with the usual absolute value.  With  $$a_K(n) = |\{\mathfrak{n} \mid N\mathfrak{n} = n\}|$$  we write $$\zeta_K(s) = \sum_{\mathfrak{n}} \frac{1}{(N\mathfrak{n})^s} =  \sum_n \frac{a_K(n)}{n^s}$$
 for the Dedekind zeta-function of $K$. Let $\mathcal{P}$ denote the set of degree one prime ideals in $\mathcal{O}_K$. Let $\Bbb{N}^{\sharp}$ denote the set of  positive integers all of whose prime factors lie over  prime ideals in $\mathcal{P}$, and $\Bbb{N}^{\flat}$ the set  of positive integers all of whose prime factors lie over  prime ideals not in $\mathcal{P}$. 
\\

\textbf{Acknowledgement:} I am very grateful to P\'eter Maga for many useful discussions and helpful suggestions on the topic of this paper. 

\section{Arithmetic in number fields}

Since $K$ is Galois, every rational prime $p$ decomposes into $r$ prime ideals of norm $p^f$ for two natural numbers $r, f$ depending on $p$ and satisfying $rf = k$. The set of primes $p$ with $f=1$, i.e.\ $p \in \mathcal{P}$, has Dirichlet density $1/k$. For an integer $n = \prod_j p_j^{\alpha_j}$ we define
\begin{equation}\label{def-n}
n^{\ast} := \prod_j p_j^{f_j \lceil \alpha_j/f_j\rceil}, \quad n^{\sharp} = \prod_{f_j = 1} p_j^{\alpha_j} \in \Bbb{N}^{\sharp}, \quad n^{\flat} = \prod_{f_j > 1} p_j^{\alpha_j} \in \Bbb{N}^{\flat},
\end{equation}
so that $n = n^{\sharp} n^{\flat}$. With this notation we obviously have
$$n^{\ast} = n \quad \text{if} \quad n \in \Bbb{N}^{\sharp},$$
and in general
\begin{equation}\label{dede-bound}
  a_K(n) \leq a_K(n^{\ast})  = \prod_j \tau_{r_j}(p^{\lceil \alpha_j/f_j\rceil})  \leq \tau_k(n)
\end{equation}
with the notation as in \eqref{def-n}. For  $n \in \Bbb{N}^{\sharp}$ we have  $\tau_{k-1}(n) = \tau_{k-1}(\mathfrak{n})$ for any ideal $\mathfrak{n}$ of norm $n$ (where the divisor function on the right is the divisor function on ideals), so that 
\begin{equation}\label{dede}
 \sum_{n \in \Bbb{N}^{\sharp}} \frac{\tau_{k-1}(n)a_K(n)}{n^s} = \zeta_K(s)^{k-1} H(s),
\end{equation}
for some Euler product $H$ that is holomorphic in $\Re s > 1/2$.

Let $\mathfrak{q}$ be an integral ideal.  Then $n \mid N\mathfrak{q}$ is equivalent to the statement that $\mathfrak{q}$ is divisible by some ideal of norm $n^{\ast}$, in particular
\begin{equation}\label{upper}
\textbf{1}_{n \mid N\mathfrak{q}} \leq \sum_{N\mathfrak{n} = n^{\ast}} \textbf{1}_{\mathfrak{n}   \mid  \mathfrak{q}}. 
\end{equation}
By inclusion-exclusion there exists a function $\mu_{n}(\mathfrak{n})$, such that
\begin{equation}\label{inclexcl}
\textbf{1}_{n \mid N\mathfrak{q}} =  \sum_{\mathfrak{n}  }  \mu_n( \mathfrak{n}) \textbf{1}_{\mathfrak{n}   \mid  \mathfrak{q}},
\end{equation}
and this function is supported on ideals $\mathfrak{n}$ satisfying   $n^{\ast} \mid N \mathfrak{n}$ and  $\mathfrak{n} \mid {\rm lcm}\{\mathfrak{q} \mid N\mathfrak{q} = n^{\ast}\}$, in particular 
\begin{equation}\label{support}
n^{\ast} \mid N\mathfrak{n} \mid n^k. 
\end{equation}
By M\"obius inversion, this function is given explicitly as
$$\mu_n(\mathfrak{q}) := \sum_{\substack{ n \mid N\mathfrak{n}\\ \mathfrak{n} \mid \mathfrak{q}}} \mu(\mathfrak{q} \mathfrak{n}^{-1}),$$
where $\mu$ denotes the usual M\"obius functions on ideals. 
 We only use this formula to conclude the trivial upper bound
\begin{equation}\label{mu-bound}
|\mu_n(\mathfrak{q})| \leq 2^{k \omega(n)},
\end{equation}
since for each prime $p \mid n$ there are at most $2^k$ squarefree ideals $\mathfrak{n}$ of $p$-power norm that can contribute to the sum non-trivially. 


For a region $\mathcal{R} \subseteq \Bbb{R}^{k-1} \setminus \{0\}$, $X > 1$ and an integral  ideal $\mathfrak{n}$ let 
$$\mathcal{A}_X(\mathfrak{n}) =    \left\{\textbf{x} \in   \Bbb{Z}^{k-1} : \textbf{x} \in \mathcal{R}_X, \, \mathfrak{n} \mid (\textbf{x}) \right\} $$
and 
\begin{equation}\label{rho}
\rho(\mathfrak{n}) = \frac{1}{(N\mathfrak{n})^{k-2} }\left|\left\{ \textbf{x} \in (\Bbb{Z}/N\mathfrak{n}\Bbb{Z})^{k-1} : \mathfrak{n} \mid (\textbf{x})\right\}\right|.
\end{equation}
This function on ideals is connected with the function $\varrho$  defined in \eqref{varrho} by
\begin{equation}\label{rhorho}
\begin{split}
\frac{\varrho(n)}{n} & = \frac{1}{n^{k-1}} \sum_{ \textbf{x} \in (\Bbb{Z}/n\Bbb{Z})^{k-1}} \textbf{1}_{n \mid f(\textbf{x})} =   \frac{1}{n^{k-1}} \sum_{ \textbf{x} \in (\Bbb{Z}/n\Bbb{Z})^{k-1}} \sum_{\mathfrak{n}  }  \mu_n( \mathfrak{n}) \textbf{1}_{\mathfrak{n}   \mid  (\textbf{x})}\\
&= \frac{1}{n^{k-1}}\sum_{\mathfrak{n}  }  \mu_n( \mathfrak{n}) \frac{n^{k-1}}{(N\mathfrak{n})^{k-1}} (N\mathfrak{n})^{k-2}\rho(\mathfrak{n}) =   \sum_{\mathfrak{n}  }  \mu_n( \mathfrak{n}) \frac{\rho(\mathfrak{n})}{N\mathfrak{n}}.
\end{split}
\end{equation}
The function $\rho$ is multiplicative in the sense that 
\begin{equation}\label{mult}
  \rho(\mathfrak{n}_1 \mathfrak{n}_2) =  \rho(\mathfrak{n}_1)\rho( \mathfrak{n}_2), \quad (N\mathfrak{n}_1, N\mathfrak{n}_2) = 1,
\end{equation}
see  \cite[Lemma 6.5]{Ma}. 
Let $N\mathfrak{n} = p^{\ell} = p^{ak + b}$ for a prime $p$, where $0 \leq b \leq k-1$ and $a \geq 0$. Since $\mathfrak{n} \mid (\textbf{x})$ implies $N\mathfrak{n} \mid N(\textbf{x})$, we have  
\begin{equation}\label{returning}
\rho(\mathfrak{n}) \leq \varrho(N\mathfrak{n}) = \varrho(p^{\ell}) = \frac{1}{p^{\ell(k-2)} }|\{ \textbf{x} \in (\Bbb{Z}/p^{\ell}\Bbb{Z})^{k-1} : p^{\ell} \mid f(\textbf{x}) \}|. 
\end{equation}
To get upper bounds on the right-hand side, it is convenient to introduce the modified counting function
$$\varrho^{\ast}(n) :=   \left| \left\{ \textbf{x} \in (\Bbb{Z}/n\Bbb{Z})^{k-1} : n \mid f(\textbf{x}), \, {\rm gcd}(x_1, \ldots, x_{k-1}, n) = 1\right\} \right|.$$
A standard application of Hensel's Lemma (\cite[Korollar 1.13]{D2}, which is also implicit in \cite[Section 10]{Ho}) shows
$$\varrho^{\ast}(p^{\alpha})  \ll p^{\alpha(k-2)}$$
with an implicit constant depending only on $f$. Returning to \eqref{returning},  we now count separately the tuples $\textbf{x} \in (\Bbb{Z}/p^{\ell}\Bbb{Z})^{k-1}$ with $p^{\alpha} \parallel {\rm gcd}(x_1, \ldots, x_{k-1})$ for $\alpha = 0, 1, \ldots, a$ and $p^{a+1} \mid {\rm gcd}(x_1, \ldots, x_{k-1})$.  Hence for fixed $\alpha$ we need to count vectors $\textbf{y} \in (\Bbb{Z}/p^{\ell - \alpha}\Bbb{Z})^{k-1}$ not all of whose components are divisible by $p$ and satisfying $p^{\ell - \alpha k}  \mid f(\textbf{y})$. This last condition is void if $\alpha > a$, so that
\begin{equation}\label{rho2}
\begin{split}
\varrho(p^{\ell})&  \leq \frac{1}{p^{\ell(k-2)}}\Bigl(\sum_{0 \leq \alpha \leq a}  p^{((\ell - \alpha)  - (\ell - \alpha k))(k-1)} \varrho^{\ast}(p^{\ell - \alpha k}) + p^{(k-1)(\ell - a - 1)} \Bigr)\\
& \ll  \frac{1}{p^{\ell(k-2)}}\Bigl(\sum_{0 \leq \alpha \leq a} p^{\alpha - 2\ell + k\ell} + p^{(k-1)(\ell - a - 1)} \Bigr)  \ll  p^a + p^{a+b+1 - k} \ll p^a = p^{[\ell/k]}.
\end{split}
\end{equation}
In connection with \eqref{returning} and \eqref{mult} we conclude in particular
\begin{equation}\label{inparticular}
  \rho(\mathfrak{n}) \ll (N\mathfrak{n})^{1/k}. 
\end{equation}
By \cite[Lemma 6.5]{Ma} we have 
\begin{equation}\label{rho3}
\rho(\mathfrak{p}) = 1  \quad \text{if $\mathfrak{p}$ is a degree 1 ideal.}
\end{equation}
Together with \eqref{rhorho}, \eqref{mu-bound} and \eqref{inparticular} we obtain
\begin{equation}\label{rho4}
\varrho(p) = k + O(p^{-1+\frac{1}{k}}), \quad p \in \mathcal{P}. 
\end{equation}
The same argument yields
  \begin{equation}\label{rho5}
 \varrho(p) \ll p \frac{(p^{\ast})^{1/k}}{p^{\ast}} \leq \frac{1}{p^{1 - 2/k}}, \quad p \not\in \mathcal{P}. 
 \end{equation}  
From \eqref{rho2}, \eqref{rho4}, \eqref{rho5} it is not hard to see that the Euler product of 
$$\zeta(s)^{1-k} \sum_n \frac{\tau_{k-1}(n)\varrho(n) }{n^s}$$
is absolutely convergent in $\Re s > \max(1/2, 2/k)$, so that in particular the product on the right-hand side of \eqref{C} is absolutely convergent.

\section{Lattice points count}\label{level}  
 
We make some choices. Fix once and for all a set $\mathfrak{C}$ of integral  ideals representing the class group of $\mathcal{O}_K$. For each integral ideal $\mathfrak{n}$, there exists a unique $\mathfrak{c} \in \mathfrak{C}$ such that $\mathfrak{n} \mathfrak{c}$ is principal. We choose, once and for all, a generator $\textbf{n} = \sum_j n_j \omega_j \in \mathcal{O}_K$ all of whose conjugates are of comparable size, i.e.\ 
\begin{equation}\label{embeddings}
|\textbf{n}^{\sigma}|  \asymp (N\mathfrak{n})^{1/k}
\end{equation}
 for all embeddings $\sigma$ of $K$ into $\Bbb{R}$ or $\Bbb{C}$, and 
\begin{equation}\label{norms}
\| \textbf{n} \|  = \|(n_1, \ldots, n_k)\| \ll   (N\mathfrak{n})^{1/k}.
\end{equation}
 This can be achieved by multiplying $\textbf{n}$ with a suitable unit of $\mathcal{O}_K$, if necessary. 

The following lemma   is a standard application of lattice reduction.
\begin{lemma} For an integral ideal $\mathfrak{n}$ there exists  $\mathbf{z} = \mathbf{z}(\mathfrak{n}) \in \Bbb{Z}^{k} \setminus \{0\}$    such that $\mathbf{z} \cdot \mathbf{n}$ has vanishing last coordinate and 
\begin{equation}\label{maynard}
| \mathcal{A}_X(\mathfrak{n})| = \frac{\rho(\mathfrak{n})}{N\mathfrak{n}} \text{{\rm vol}}(\mathcal{R}_X)+ O\left(1 + \frac{X^{k-2}}{\| \mathbf{z}(\mathfrak{n}) \|^{k-2} (N\mathfrak{n})^{(k-2)/k}} \right)
 \end{equation}
 as well as 
 \begin{equation}\label{maynard1}
| \mathcal{A}_X(\mathfrak{n}) | \ll   \frac{\rho(\mathfrak{n})}{N\mathfrak{n}} X^{k-1}+  \sum_{j=1}^{k-2} \frac{X^j}{\| \mathbf{z}(\mathfrak{n}) \|^j (N\mathfrak{n})^{j/k}} .
 \end{equation}
 \end{lemma}
 
\textbf{Proof.} The formula  \eqref{maynard} is \cite[(6.2) and subsequent display]{Ma} with $(n, k, Q, \mathfrak{d}) \mapsto (k, 1, 1, \mathfrak{n})$. The upper bound \eqref{maynard1} is a small variation, based on the fact that $0 \not\in \mathcal{R}$. For convenience we provide the details. First we observe that
$| \mathcal{A}_X(\mathfrak{n}) | \leq | \mathcal{A}_{X\cdot N\mathfrak{c}}((\textbf{n})) |$. The right-hand side counts integral vectors $\textbf{b} \in \Bbb{Z}^{k}$ satisfying $\textbf{n} \cdot \textbf{b} \in \mathcal{R}_{X \cdot N \mathfrak{c}}$. In particular, $\textbf{b}$ lies in a rank $k-1$ lattice
$$\Lambda(\textbf{n}) = \{\textbf{b} \in \Bbb{R}^{k} \mid \textbf{n} \cdot \textbf{b} \text{ has vanishing last coordinate}\}.$$
Since $\| \textbf{n} \cdot \textbf{b}  \| \ll X$, we have $|\textbf{n}^{\sigma}|  |\textbf{b}^{\sigma}| = |(\textbf{n} \cdot \textbf{b})^{\sigma}| \ll X$ for all embeddings $\sigma$, and by \eqref{embeddings} this implies $| \textbf{b}^{\sigma} | \ll X (N \mathfrak{n})^{-1/k}$ for all embeddings $\sigma$,  so $\| \textbf{b} \| \ll X (N\mathfrak{n})^{-1/k}$. 

There exists a $\Bbb{Z}$-basis $\textbf{z}_1(\textbf{n}), \ldots, \textbf{z}_{k-1}(\textbf{n})$ for $\Lambda(\textbf{n})$, such that   
 if $\textbf{b} = \sum_j b_j \textbf{z}_j(\textbf{n}) \in \Lambda(\textbf{n})$, then  $b_j \ll \| \textbf{b} \| /\| \textbf{z}_j(\textbf{n})\|$, see e.g.\ \cite[Lemma 5]{Da}. We order it such that 
$\| \textbf{z}_1(\textbf{n}) \| \leq  \ldots \leq \| \textbf{z}_{k-1}(\textbf{n})\|$. We have the essentially trivial inequality $\prod_j \| \textbf{z}_j(\textbf{n})\|  \geq \det \Lambda(\textbf{n})$ (Hadamard's inequality), and we recall that $0 \not\in \mathcal{R}$.  
Therefore we have  a total of 
$$\ll \prod_{j=1}^{k-1} \left(1 +  \frac{X }{(N \mathfrak{n})^{1/k} \| \textbf{z}_j(\textbf{n})\|}\right) - 1 \ll  \frac{X^{k-1}}{\det\Lambda(\textbf{n})} + \sum_{j=1}^{k-2} \frac{X^j}{\| \mathbf{z}_1(\mathfrak{n}) \|^j (N\mathfrak{n})^{j/k}} $$
choices for the $(k-1)$-tuple $(b_1, \ldots, b_{k-1}) \not= 0$ if $\textbf{b} = \sum_j b_j \textbf{z}_j(\textbf{n}) \in \Lambda(\textbf{n})$  and $\| \textbf{b} \| \ll X (N\mathfrak{n})^{-1/k}$. 
Letting $X \rightarrow \infty$ and using the definition \eqref{rho}, one confirms 
$$\frac{X^{k-1}}{\det\Lambda(\textbf{n})}  = \frac{\rho((\textbf{n}))}{N(\textbf{n})} X^{k-1} \ll  \frac{\rho(\mathfrak{n})}{N\mathfrak{n}} X^{k-1},$$
and the claim follows. \\

\textbf{Remark:} The key point of \eqref{maynard1} is  that the right-hand side has no contribution $O(1)$ for $j=0$. This feature is needed in the passage from \eqref{e21a} to \eqref{e21b} below.

 \section{Upper bounds}

 An important input for the proof is the following result which is essentially due to Wolke. 
\begin{lemma}\label{wolke} Let $c, C >0$ be constants. Let $F, \kappa$ be  multiplicative functions such that $0 \leq F(p^{\alpha})  \ll  \alpha^{C}$,   
\begin{equation}\label{new}
0 \leq \kappa(p^{\alpha})   \ll p^{[\alpha/3]}
\end{equation}
 and  $\kappa(p) <p$ for all primes $p$. For a sequence $a_1, a_2, \ldots$   of natural numbers, $x \geq 1$, $d \in \Bbb{N}$ define
$$R(x, d) := \sum_{\substack{n \leq x\\ d \mid a_n}} 1 - \frac{\kappa(d)}{d} x,$$
and suppose that
\begin{equation}\label{suppose}
\sum_{d \leq x^c} |R(x, d)| \ll x^{1-c}. 
\end{equation}
Then
$$\sum_{n \leq x} F(a_n) \ll_{C, c, f, \kappa} x \exp\Bigl( \sum_{p \leq x} \frac{\kappa(p)}{p}(F(p) - 1)\Bigr).$$  
 \end{lemma}

\textbf{Proof.} This is \cite[Satz 1]{Wo} except that \eqref{suppose} is slightly stronger than \cite[($A_3$), ($A_4$)]{Wo}, but    \eqref{new} is a weaker assumption than the analogous bound  $0 \leq \kappa(p^{\alpha}) \ll \alpha^C$ in \cite[($A_2$)]{Wo}. 
  A careful inspection of the proof shows that the same argument works verbatim under the present assumption \eqref{new} (and in fact even weaker estimates are possible). In the following we list the places in the proof of \cite[Satz 1]{Wo} where bounds for $\kappa(p^{\alpha})$ with  $\alpha \geq 2$ are needed.  
\begin{itemize}
\item Lemma 3 in \cite{Wo} requires $g(p^l) \ll l^{O(1)}$ for a multiplicative arithmetic function $g$, but what is really needed to make the first display in the proof valid is $g(p) \ll 1$ and 
\begin{equation}\label{really}
\sum_{\alpha \geq 2} \sum_{p \text{ prime}} \frac{g(p^{\alpha})}{p^{\alpha(1-\eta)}} \ll 1
\end{equation}
for a suitably small $\eta > 0$. We need to apply this lemma with the assumption \eqref{really} twice in the following. 
\item For \cite[(3.3)]{Wo} it is needed that\footnote{of course, even the upper bound $z^{-\varepsilon}$ for any $\varepsilon > 0$ would suffice}
$$\sum_{r \leq \log z}\sum_{z^{1/(2(r+1))} \leq p \leq z^{1/(2r)}} \frac{\kappa(p^{r+1})}{p^{r+1}} \ll z^{-1/9}$$
for $z \geq 2$, which follows easily from \eqref{new}. 
\item For \cite[(3.4)]{Wo} it suffices that $\kappa(n) \ll n^{1/2}$ which is weaker than \eqref{new}. 
\item The bound \cite[(3.10)]{Wo} applies \cite[Lemma 3]{Wo} to the arithmetic function
$$g : b \mapsto F(b) \prod_{p^{\ell} \parallel b} \frac{p}{p - \kappa(p)} \cdot \left(\kappa(p^{\ell}) - \frac{\kappa(p^{\ell+1})}{p}\right),$$
which satisfies \eqref{really} by \eqref{new} and our assumption on $F$. 
 \item For the first sum in \cite[(3.14)]{Wo}, any polynomial bound on $\kappa(n)$ suffices in view of our assumption \eqref{suppose}, while for the second sum in \cite[(3.14)]{Wo} we 
 apply \cite[Lemma 3]{Wo} with $g = \kappa$, which satisfies \eqref{really} by \eqref{new}. \\
\end{itemize}

As a very special case we obtain for the sequence $a_n = n$ the upper bound (which can easily be proved in many other ways) 
\begin{equation}\label{divisor} 
\sum_{n \leq x} \tau_{\ell}(n)^{\beta} \ll_{\beta, \ell}  x (\log x)^{\ell^{\beta} - 1}
\end{equation}
for $\beta \geq 0$, $\ell \in \Bbb{N}$.  A more advanced consequence is the following.
\begin{cor}\label{wolke1} Let $F$ be as in the previous lemma and suppose that $F(p) = k$ for all primes $p \in \mathcal{P}$ and $F(p) = 0$ for $p \not\in \mathcal{P}$.   Let $V \geq 2$. Then
$$\sum_{\substack{0 \not= \mathbf{v} \in \Bbb{Z}^{k-1} \\ \| \mathbf{v} \| \leq V}} F(N(\mathbf{v})) \ll V^{k-1} (\log V)^{k-1}.$$
\end{cor}

\textbf{Proof:} We order the natural numbers $N(\textbf{v})$, $0 \not= \textbf{v} \in \Bbb{Z}^{k-1}$, by size and call this the sequence $a_1, a_2, \ldots$. For $d \in \Bbb{N}$ we conclude from  \eqref{inclexcl},  \eqref{maynard}, \eqref{mu-bound} and \eqref{support} that
$$\sum_{\substack{0 \not= \textbf{v} \in \Bbb{Z}^{k-1}\\  \| \textbf{v} \| \leq V \\ d \mid N(\textbf{v})}} 1 = \sum_{\mathfrak{n}} \mu_d(\mathfrak{n}) \sum_{\substack{0 \not= \textbf{v} \in \Bbb{Z}^{k-1}\\  \| \textbf{v} \| \leq V \\ \mathfrak{n} \mid (\textbf{v})}} 1 = \sum_{\mathfrak{n}} \mu_d(\mathfrak{n})\Bigl( \frac{\rho(\mathfrak{n})}{N\mathfrak{n}} (2V)^{k-1} + O(V^{k-2})  \Bigr) = \frac{\varrho(d)}{d} (2V)^{k-1} +O_{\varepsilon}(d^{\varepsilon} V^{k-2})$$
for any $\varepsilon > 0$, from which we obtain easily that the sequence $a_n$ satisfies \eqref{suppose} with $\kappa = \varrho$. The bound \eqref{new} holds by \eqref{rho2}, and $\varrho(p) < p$ is clear from the fact that the incomplete norm form $f(\textbf{x})$ has no fixed divisor, for instance $f((1, 0, \ldots, 0)) = 1$. Hence Lemma \ref{wolke} is applicable, and we conclude the desired bound from
$$\sum_{p \leq x} \frac{\varrho(p)}{p}(F(p) - 1)  \ll (k-1)\log \log x$$
by \eqref{rho4} and the fact that $\mathcal{P}$ has Dirichlet density $1/k$. 
 
\section{Decomposition of divisor functions} 
 
In this section we decompose the divisor function $\tau_k$ following ideas of \cite{FI}, which is a somewhat sophisticated generalization of Dirichlet's hyperbola method. 
Let $y \geq 1$.  Then  we have by inclusion-exclusion
$$\tau_k(n) =  \sum_{\substack{n_1 \cdots n_k = n\\  n_1, \ldots, n_k \leq y}} 1 + \sum_{j=1}^{k} (-1)^{j-1}  \left(\begin{array}{l} k\\ j \end{array} \right) \sum_{\substack{n_1 \cdots n_k = n \\   n_1, \ldots, n_j > y}} 1.$$
Now let $X \geq 1$ be such that   $y \geq 2X $. If $X^k \leq n \leq 2X^k$,  then the term corresponding to $j = k$ is empty, and the condition $n_1 >y$ implies automatically $n_2 \cdots n_k \leq 2X^k y^{-1}$. In the first term we may assume, by symmetry, that $n_1$ is the largest variable. Choosing $y= X \Delta$ with $2 \leq \Delta \leq X^{1/100}$, say, we obtain
$$\tau_k(n) = \sum_{j=1}^{k-1} (-1)^{j-1} \left(\begin{array}{l} k\\ j \end{array} \right) \sum_{\substack{n_1 \cdots n_k = n \\ n_2 \cdots n_k \leq 2X^{k-1} \Delta^{-1}\\ n_1, \ldots, n_j > X \Delta}} 1 + O\Bigl(\sum_{\substack{n_1 \cdots n_k = n\\ X^{k-1} \Delta^{-1}  \leq n_2 \cdots  n_k  \leq 2X^{k-1} 
}} 1 \Bigr).$$
We can drop the condition $n_1 > X\Delta$ in the main terms, for if $n_1 \leq X\Delta$, then $n_2 \cdots n_k \geq X^{k-1} \Delta^{-1}$, and  this contribution can be absorbed in the error term. 

This gives 
\begin{equation}\label{mainterm}
\mathcal{M}(\mathcal{R}_X) = \sum_{j=1}^{k-1} (-1)^{j-1} \left(\begin{array}{l} k\\ j \end{array} \right) M_j + O(E) ,
\end{equation}
where
$$M_j = \sum_{\substack{  n_2 \cdots n_k \leq 2X^{k-1} \Delta^{-1}\\ n_2, \ldots, n_j > X \Delta}} \sum_{\substack{\textbf{x} \in \mathcal{R}_X \cap \Bbb{Z}^{k-1} \\  n_2 \cdots n_k \mid f(\textbf{x})}} 1, \quad\quad  E  = \sum_{X^{k-1} \Delta^{-1}  \leq n  \leq 2X^{k-1}}  \tau_{k-1}(n) \sum_{\substack{\textbf{x} \in \mathcal{R}_X \cap \Bbb{Z}^{k-1} \\   n  \mid f(\textbf{x})}} 1. 
$$

For later purposes it is convenient to smooth the sums $M_j$. Let $W^{+}$ be a fixed non-negative smooth function that is 1 on $[0, 2]$ and 0 on $[3, \infty)$, and let $W^{-}$ be a fixed non-negative smooth function that is 1 on $[0, 1]$ and 0 on $[2, \infty)$. Similarly let $V^{+}$ be a fixed non-negative smooth function that is 1 on $[0, 1/2]$ and 0 on $[1, \infty)$, and let $V^{-}$ be a fixed non-negative smooth function that is 1 on $[0, 1]$ and 0 on $[2, \infty)$. The Mellin transforms of $V^{\pm}$ and $W^{\pm}$ have simple poles with residue 1 at $s = 0$ and are rapidly decaying on vertical lines. 

We clearly have
\begin{equation}\label{inequality}
M_j^{-} \leq M_j \leq M_j^{+},
\end{equation}
where
\begin{equation}\label{mpm}
M_j^{\pm} = \sum_{   n_2,  \ldots, n_k}  W^{\pm}\left(\frac{n_2 \cdots n_k}{X^{k-1} \Delta^{-1}}\right)   \prod_{i=2}^j V^{\pm} \left(\frac{X\Delta}{n_i}\right) \sum_{\substack{\textbf{x} \in \mathcal{R}_X \cap \Bbb{Z}^{k-1} \\  n_2 \cdots n_k \mid f(\textbf{x}) }} 1.
\end{equation}
As usual, an empty product is interpreted as 1. 

\section{Error terms I}\label{sec6}

This section is devoted to the estimation of the error term $E$ in \eqref{mainterm}. The final bound is \eqref{finalE} below. 

By \eqref{upper} and \eqref{maynard1} we have 
\begin{displaymath}
\begin{split}
E & \leq \sum_{  X^{k-1} \Delta^{-1} \leq   n  \leq 2X^{k-1} 
} \tau_{k-1}(n)   \sum_{N\mathfrak{n} =n^{\ast}} |\mathcal{A}_X(\mathfrak{n} )|   \ll   E_1 + E_2 ,
\end{split}
\end{displaymath}
where
\begin{displaymath}
\begin{split}
E_1 & =  X^{k-1} \sum_{  X^{k-1} \Delta^{-1} \leq   n  \leq 2X^{k-1}} \tau_{k-1}(n) \sum_{N\mathfrak{n} =n^{\ast}} \frac{\rho(\mathfrak{n})}{N\mathfrak{n}} ,\\
E_2 &= \sum_{j=1}^{k-2} X^j \sum_{n  \leq 2X^{k-1}} \tau_{k-1}(n) \sum_{N\mathfrak{n} =n^{\ast}} \frac{1}{\| \textbf{z}(\mathfrak{n})\|^j (N\mathfrak{n})^{j/k}}.
\end{split}
\end{displaymath}
We start with the estimation of $E_1$. We decompose uniquely $n = n_1 n_2 m$ with $n_1n_2 \in \Bbb{N}^{\sharp}$, $n_1$ squarefree, $n_2$ squarefull, $n_1, n_2$ coprime , $m \in \Bbb{N}^{\flat}$. This notation in effect, we infer from \eqref{mult}, \eqref{inparticular}, \eqref{rho3} and \eqref{dede-bound} that 
\begin{equation}\label{E1}
\begin{split}
E_1 & \ll X^{k-1} \sum_{\substack{X^{k-1} \Delta^{-1} \leq n_1n_2 m \leq 2 X^{k-1} }}   \frac{\tau_{k-1}(n_2m)a_K(n_2m^{\ast})}{(n_2m^{\ast})^{1-1/k}}\cdot \frac{\tau_{k-1}(n_1) a_K(n_1)}{n_1}\\
& \ll X^{k-1} (\log X)^{k-2} \log \Delta \sum_{n_2,m}   \frac{\tau_{k-1}(n_2m)\tau_k(n_2m^{\ast})}{(n_2m^{\ast})^{1-1/k}} \ll X^{k-1} (\log X)^{k-2} \log \Delta. 
\end{split}
\end{equation}

The estimation of $E_2$ is the most delicate part of the argument, since we have not even a logarithm to spare.  To estimate $E_2$ we let $B \in \Bbb{N}$ be a very large constant, 
$$  \alpha(n) := \min\left(\tau_{k}(n)^2 \tau_B(n), \exp((\log\log X)^2)\right),$$
and 
$$Z_0(n) := X^{1/k} (\log X)^{-\frac{1}{k-1}} \left(\frac{n}{n^{\ast}}\right)^{1/k}\alpha(n^{\flat}).$$
This rather artificial definition is carefully designed   and takes care in particular of the contribution of prime ideals of degree $> 1$. 

We split $E_2$ into two subsums $E_{21}$ and $E_{22}$ according to whether $\| \textbf{z}(\mathfrak{n})\| \geq Z_0(n)$ or $\| \textbf{z}(\mathfrak{n})\|< Z_0(n)$, respectively. By \eqref{dede-bound} we have
\begin{equation}\label{e21a}
\begin{split}
E_{21}& \leq \sum_{j=1}^{k-2} X^j \sum_{n  \leq 2X^{k-1}} \frac{ \tau_{k-1}(n)  a_K(n^{\ast})}{ Z_0(n)^j (n^{\ast})^{j/k}}\\
& \leq \sum_{j=1}^{k-2} X^{j(1 - \frac{1}{k})} (\log X)^{\frac{j}{k-1}}  \sum_{  m \in \Bbb{N}^{\flat}} \frac{\tau_{k-1}(m) \tau_k(m)}{m^{j/k} \alpha(m)^j} \sum_{\substack{n  \leq 2X^{k-1}/m\\ n \in \Bbb{N}^{\sharp}}}  \frac{\tau_{k-1}(n)  a_K(n)}{n^{j/k}}.  
\end{split}
\end{equation}
We estimate the inner sum using \eqref{dede}, and it follows easily that
\begin{equation}\label{e21b}
E_{21} \ll \sum_{j=1}^{k-2} X^{j(1-\frac{1}{k})} (\log X)^{\frac{j}{k-1}} \sum_{  m  \leq 2X^{k-1}} \frac{\tau_{k-1}(m) \tau_k(m)}{m^{j/k} \alpha(m)}\Bigl( \frac{X^{k-1}}{m}\Bigr)^{ 1 - \frac{j}{k}} (\log X)^{k-2}.
\end{equation}
Using \eqref{divisor}, we estimate the $m$-sum by
\begin{displaymath}
\begin{split}
\sum_{m \leq 2X^{k-1}}  \frac{\tau_{k-1}(m) \tau_k(m)}{m \,\alpha(m)}  \leq \sum_{ m  \leq 2X^{k-1}}  \frac{ 1}{m \,\tau_B(m)}   + \sum_{ m \leq 2X^{k-1}  } \frac{\tau_k(m)^2}{m \exp((\log\log X)^2)} \ll_B (\log X)^{1/B}, 
\end{split}
\end{displaymath}
so that
\begin{equation}\label{E21}
\begin{split}
E_{21} & \ll_B X^{k-1} (\log X)^{k-1 - \frac{1}{k-1} + \frac{1}{B}}.
\end{split}
\end{equation}

We now turn towards the estimation of 
$$E_{22} = \sum_{j=1}^{k-2} X^j\sum_{ n  \leq 2X^{k-1}}  \tau_{k-1}(n)  \sum_{\substack{ N\mathfrak{n} = n^{\ast}\\  \| \textbf{z} ( \mathfrak{n}) \| < Z_0(n)}} \frac{1}{\| \textbf{z} (  \mathfrak{n}) \|^j (N\mathfrak{n})^{j/k}}.$$
Let $\textbf{n} \in \mathcal{O}_K$ be as in the beginning of Section \ref{level}. The idea is now to glue together $\textbf{z} = \textbf{z}(\mathfrak{n})$ and $\textbf{n}$ and to consider the  non-zero integral vector  $\textbf{v} = \textbf{z}\cdot  \textbf{n}$, which by definition has vanishing last coordinate. By \eqref{norms} we have
\begin{displaymath}
\begin{split}
\| \textbf{v}  \|&  \ll Z_0(n) \cdot  (n^{\ast})^{1/k} \leq  
2^{1/k} X(\log X)^{-1/(k-1)} \alpha( n^{\flat})  . 
%
\end{split}
\end{displaymath}
Let $E_{221}$ denote the contribution of those $\mathfrak{n}$ where $\| \textbf{v} \| \leq V_0 := X (\log X)^{-1/(k-1)}$, so that
\begin{displaymath}
\begin{split}
E_{221}&  \leq \sum_{j=1}^{k-2} X^j \sum_{\substack{0 \not= \textbf{v} \in \Bbb{Z}^{k-1}\\ \| \textbf{v} \| \leq V_0}}  \frac{R((\textbf{v}))}{\| \textbf{v} \|^j}\end{split}
\end{displaymath}
with
$$R(\mathfrak{q}) := \sum_{\mathfrak{n} \mid \mathfrak{q}} \sum_{n : n^{\ast} = N\mathfrak{n}} \tau_{k-1}(n).
$$
Using the simple bound 
  $|\{\mathfrak{n} \mid \mathfrak{q}\} | \leq \tau(N\mathfrak{q})^k$, we have
  $$R(\mathfrak{q}) \leq  \tau(N\mathfrak{q})^{k+1} \tau_{k-1}(N\mathfrak{q})$$
in general, and 
$$R(\mathfrak{q}) = k, \quad N\mathfrak{q} = p \in \mathcal{P},$$
while the case $N\mathfrak{q} =  p \not\in \mathcal{P}$ cannot occur. Clearly $R$ is multiplicative in the sense that $R(\mathfrak{q}_1\mathfrak{q}_2) = R(\mathfrak{q}_1)R(\mathfrak{q}_2)$ if $(N\mathfrak{q}_1, N\mathfrak{q}_2) = 1$. 
 Let $T : \Bbb{N} \rightarrow \Bbb{N}$ be the multiplicative function defined by
$$T(p) = \begin{cases} k, & p \in \mathcal{P},\\0, & p \not\in \mathcal{P}, \end{cases} \quad\quad  T(p^{\alpha}) =  \tau(p^\alpha)^{k+1} \tau_{k-1}(p^{\alpha})$$
for primes $p$ and $\alpha \in \{2, 3, \ldots\}$. Then by Corollary \ref{wolke1} we obtain
\begin{equation}\label{E221}
\begin{split}
 E_{221}&\leq \sum_{j=1}^{k-2} X^j \sum_{2^{\nu} \leq V_0} \frac{1}{2^{j \nu}} \sum_{\substack{ \textbf{v} \in \Bbb{Z}^{k-1}\\ 2^{ \nu} \leq \| \textbf{v} \| \leq 2^{\nu+1}    }} T(N(\textbf{v}))\\
 &  \ll  \sum_{j=1}^{k-2} X^j \sum_{2^{\nu} \leq V_0} \frac{1}{2^{j \nu}}  2^{\nu(k-1)} (\log 2^{\nu})^{k-1} \ll  \sum_{j=1}^{k-2} X^j  V_0^{k-1 - j} (\log X)^{k-1}\\
 & \ll X^{k-1} (\log X)^{k-1 - \frac{1}{k-1}}.
 \end{split}
\end{equation}


Denoting by $E_{222}$ the remaining contribution with $V_0 \leq \| \textbf{v} \| \ll V_1 := X \exp((\log \log X)^2)$, we have by a version of Rankin's trick applied to the condition $\| \textbf{v} \| = \| \textbf{z} \cdot \textbf{n}\| \ll  V_0 \alpha(n^{\flat}) $ that
\begin{displaymath}
\begin{split}
E_{222} & \ll  \sum_{j=1}^{k-2} X^j \sum_{V_0 \leq 2^{\nu} \ll V_1} \sum_{ n  \leq 2X^{k-1}}  \tau_{k-1}(n)  \sum_{\substack{ N\mathfrak{n} = n^{\ast}\\  2^{\nu} \leq \| \textbf{z} \cdot \textbf{n}\| \leq 2^{\nu+1}\\ \| \textbf{z} \cdot \textbf{n}\| \ll  V_0 \alpha(n^{\flat})  }} \frac{1}{\| \textbf{z} (  \mathfrak{n}) \|^j (N\mathfrak{n})^{j/k}}\\
& \ll \sum_{j=1}^{k-2} X^j \sum_{V_0 \leq 2^{\nu} \ll V_1} \sum_{\substack{\textbf{v} \in \Bbb{Z}^{k-1} \\ 2^{\nu}\leq  \| \textbf{v} \| \leq 2^{\nu+1}}}  \frac{1}{\| \textbf{v} \|^j}   \sum_{\mathfrak{n} \mid (\textbf{v})} \sum_{n : n^{\ast} = N\mathfrak{n}} \tau_{k-1}(n)   \Bigl(\frac{V_0 \alpha(n^{\flat}) }{\| \textbf{v} \| }\Bigr)^{k-1-j}.
\end{split}
\end{displaymath}
Recalling the definition of $V_0 = X (\log X)^{-1/(k-1)}$, we obtain
\begin{displaymath}
\begin{split}
E_{222}  & \ll X^{k-1} (\log X)^{-\frac{1}{k-1}} \Bigl(\log\frac{V_1}{V_0}\Bigr)  \max_{V_0 \leq V \ll V_1}  \frac{1}{V^{k-1}} 
 \sum_{\substack{0 \not= \textbf{v} \in \Bbb{Z}^{k-1} \\   \| \textbf{v} \| \ll V }}  \Bigl( \sum_{\mathfrak{n} \mid (\textbf{v})} \sum_{n : n^{\ast} = N\mathfrak{n}} \tau_{k-1}(n) \alpha(n^{\flat})^{k-2}\Bigr).  
\end{split}
\end{displaymath}
The last parenthesis can be estimated by $\tilde{T}(N(\textbf{v}))$, where $\tilde{T} : \Bbb{N} \rightarrow \Bbb{N}$ is the multiplicative function defined by
$$\tilde{T}(p) = T(p), \quad \tilde{T}(p^{\alpha}) = T(p^{\alpha}) \left(\tau_k(p^{\alpha})^2 \tau_B(p^{\alpha})\right)^{k-2}$$
for primes $p$ and $\alpha \in \{2, 3, \ldots\}$. 
Again by Corollary \ref{wolke1} we obtain
\begin{equation}\label{e222}
  E_{222} \ll_B  X^{k-1} (\log X)^{-\frac{1}{k-1}} (\log\log X)^2  (\log X)^{k-1}.
\end{equation}
We summarize 
 \eqref{E1}, \eqref{E21} -- \eqref{e222}  by stating that 
\begin{equation}\label{finalE}
  E \ll_{\varepsilon} X^{k-1}\left( (\log X)^{k-2} \log \Delta + (\log X)^{k - 1 - \frac{1}{k-1} + \varepsilon}\right)
\end{equation}
where $\varepsilon = 1/B$.

 \section{Error terms II}
 
  Next we investigate the main terms $M_j^{\pm}$, defined in \eqref{mpm}.   
   By \eqref{inclexcl} 
   we have
\begin{displaymath}
\begin{split}
M_j^{\pm}& = \sum_{  n_2, \ldots, n_k } W^{\pm}\left(\frac{n_2 \cdots n_k}{X^{k-1} \Delta^{-1}}\right)   \prod_{i=2}^j V^{\pm} \left(\frac{X\Delta}{n_i}\right) \sum_{\mathfrak{n}}  \mu_{n_2\cdots n_k} (\mathfrak{n}) |\mathcal{A}_X(\mathfrak{n})|.\\
\end{split}
\end{displaymath}
 We use \eqref{maynard} for the evaluation of $|\mathcal{A}_X(\mathfrak{n})|$. The  aim of this section is  to handle the two error terms whose total contribution we call $F_1$ and $F_2$.  By \eqref{support},  \eqref{mu-bound} and \eqref{dede-bound}, the first error term contributes at most
\begin{equation}\label{E-first} 
\begin{split}
F_1 &  \ll \sum_{  n \leq 3X^{k-1} \Delta^{-1} }  \tau_{k-1}(n) \sum_{   N \mathfrak{n}   \mid n^k  } 2^{k \omega(n)} \\
&\leq \sum_{  n \leq 3X^{k-1} \Delta^{-1} }  \tau_{k-1}(n)   2^{k \omega(n)} \sum_{j=0}^k \tau_k(n^j)
 = \frac{X^{k-1} }{\Delta} (\log X)^{O(1)}. 
  \end{split}
 \end{equation}
The second error term contributes at most
\begin{equation*}
F_2 \ll  \sum_{n \leq 3X^{k-1} \Delta^{-1} }  \tau_{k-1}(n)  \sum_{\mathfrak{n}}  |\mu_{n} (\mathfrak{n})| \frac{X^{k-2}}{\| \textbf{z}(\mathfrak{n})\|^{k-2}(N\mathfrak{n})^{(k-2)/k}} .  
\end{equation*}
Again we distinguish two cases depending on whether $\|\textbf{z}(\mathfrak{n})\|$ is big or not, but the present situation is more relaxed and a slightly simpler argument than in the previous section suffices.  Let
$$Z_1(\mathfrak{n}):= X (N\mathfrak{n})^{-\frac{1}{k}} \Delta^{-\frac{1}{k-1}}.$$
The portion $F_{21}$, say, with $\|\textbf{z}(\mathfrak{n})\| \geq Z_1(\mathfrak{n})$ can be estimated in the same way as in \eqref{E-first} by
\begin{equation}\label{F21}
F_{21} \leq  \Delta^{1 - \frac{1}{k-1}}\sum_{n \leq 3X^{k-1} \Delta^{-1} }  \tau_{k-1}(n) \sum_{   N \mathfrak{n}   \mid n^k  } 2^{k \omega(n)}  = \frac{X^{k-1} }{\Delta^{1/(k-1)}} (\log X)^{O(1)}.
\end{equation}
For the portion $F_{22}$ with $\|\textbf{z}(\mathfrak{n})\| < Z_1(\mathfrak{n})$ we define again $\textbf{n} \in \mathcal{O}_K$ as in Section \ref{level} and  create the new vector $\textbf{v} := \textbf{z}(\mathfrak{n})\cdot  \textbf{n} $ with vanishing last coordinate and 
$\| \textbf{v}\| \ll X\Delta^{-1/(k-1)}$. Recalling \eqref{support} and \eqref{mu-bound} we obtain
\begin{displaymath}
\begin{split}
F_{22} & \ll X^{k-2} \sum_{\substack{\textbf{v} \in \Bbb{Z}^{k-1}\\ \| \textbf{v} \| \leq X\Delta^{-1/(k-1)}}} \frac{1}{\| \textbf{v} \|^{k-2}} \sum_{\mathfrak{n} \mid (\textbf{v})} \sum_{ n  \mid N\mathfrak{n}  } \tau_{k-1}(n) 2^{k\omega(n)} \\
&\leq  X^{k-2} \sum_{2^{\nu} \leq X \Delta^{-1/(k-1)}} \frac{1}{2^{\nu(k-2)}} \sum_{\substack{\textbf{v} \in \Bbb{Z}^{k-1}\\ 2^{\nu} \leq \| \textbf{v} \| \leq 2^{\nu+1}}}  (\tau^2_2\tau_{k-1}\tau_k)(N(\textbf{v}))  . 
\end{split}
\end{displaymath}
 By  Corollary \ref{wolke1}   we obtain
  \begin{equation}\label{F22}
  F_{22} \ll \frac{X^{k-1} }{\Delta^{1/(k-1)}} (\log X)^{O(1)}. 
 \end{equation}

 \section{Multiple $L$-functions} 
 
Collecting the   error terms   \eqref{E-first} -- \eqref{F22}, we see that
\begin{equation}\label{error-F}
M_j^{\pm}  =  \text{vol}(\mathcal{R}_X)   \tilde{M}_j^{\pm} + O\left(X^{k-1} (\log X)^{O(1)} \Delta^{-1/(k-1)} \right)
\end{equation}
where 
\begin{displaymath}
\begin{split}
 \tilde{M}^{\pm}_j &  :=   \sum_{  n_2, \ldots, n_k } W^{\pm}\left(\frac{n_2 \cdots n_k}{X^{k-1} \Delta^{-1}}\right)   \prod_{i=2}^j V^{\pm} \left(\frac{X\Delta}{n_i}\right)   \sum_{\mathfrak{n}}  \mu_{n_2 \cdots n_k}(\mathfrak{n}) \frac{\rho(\mathfrak{n})}{N\mathfrak{n} }\\
 & =  \sum_{  n_2, \ldots, n_k } W^{\pm}\left(\frac{n_2 \cdots n_k}{X^{k-1} \Delta^{-1}}\right)   \prod_{i=2}^j V^{\pm} \left(\frac{X\Delta}{n_i}\right) \frac{\varrho(n_2 \cdots n_k ) }{n_2\cdots n_k }.
  \end{split}
 \end{displaymath}
 The last equality follows from \eqref{rhorho}.  Let
$$L(s_2, \ldots, s_k) := \sum_{n_2, \ldots, n_k} \frac{\varrho(n_2 \cdots n_k)}{n_2^{s_2} \cdots n_k^{s_k}}.$$
It follows easily from \eqref{rho2}, \eqref{rho4} and \eqref{rho5}   
that the function $L(s_2, \ldots, s_k) $ is absolutely convergent in $\Re s_2, \ldots, \Re s_k > 1$, and that the Euler product
$$L(s_2, \ldots, s_k)  \prod_{j=2}^k \zeta_K(s_j)^{-1} = \prod_{p}\Bigl(1 +O\Bigl(   \sum_{j=2}^k \frac{p^{2/k}}{p^{2\Re s_j}} + \frac{p^{2/k}}{p^{\Re s_j + 1}} \Bigr) \Bigr)$$
is absolutely convergent in $\Re s_2, \ldots, \Re s_k > 1/2  +1/k $. In particular, $L(s_2, \ldots, s_k)$ can be continued meromorphically to this region   with polynomial growth on vertical lines and  polar lines at $s_i = 1$, $i = 2, \ldots, k$. We have
$$ \underset{s_2 = \ldots =  s_{k} = 1}{\text{res}} L(s_2, \ldots, s_{k} )   = C$$
with $C$ as in \eqref{C}.  
By Mellin inversion we have 
\begin{displaymath}
\begin{split}
\tilde{M}^{\pm}_j  =   \int _{(\varepsilon)}\cdots &\int_{(\varepsilon)} \int_{(2\varepsilon)}  \widehat{W}^{\pm}(s_1) \prod_{i=2}^j \widehat{V}^{\pm}(s_i) \frac{(X^{k-1} \Delta^{-1})^{s_1}}{  (X\Delta)^{s_2 + \ldots + s_j}}\\
& \times   L(1 + s_1 - s_2, \ldots, 1 + s_1 - s_j, 1 + s_1, \ldots, 1 + s_1) \frac{{\rm d}s_1 \, {\rm d}s_2 \cdots {\rm d}s_j}{(2\pi i)^j}. 
  \end{split}
  \end{displaymath}
  We shift the $s_2, \ldots, s_j$ contours  to the  \emph{right} to $\Re s_i = 1/8$, say. For each integration variable we pick up a residue at $s_i = s_1$ and obtain a remaining integral that we bound trivially, recalling that $\widehat{W}^{\pm}$ and $\widehat{V}^{\pm}$ are rapidly decaying on vertical lines. In this way we obtain 
 \begin{displaymath}
\begin{split}
\tilde{M}^{\pm}_j  =     \int_{(2\varepsilon)}  &\widehat{W}^{\pm}(s_1)    \widehat{V}^{\pm}(s_1)^{j-1} \frac{(X^{k-1} \Delta^{-1})^{s_1}}{  (X\Delta)^{(j-1)s_1}}\\
&  \underset{u_2 = \ldots =  u_{j} = 1}{\text{res}} L(u_2, \ldots, u_{j}, 1 + s_1, \ldots, 1 + s_1) \frac{{\rm d}s_1 }{(2\pi i)^j} + O\Bigl(X^{-\frac{1}{8} +  \varepsilon}\Bigr)
  \end{split}
  \end{displaymath} 
for any $\varepsilon > 0$. Next we shift the $s_1$-contour to the left past $\Re s = 0$.  Since $\widehat{W}^{\pm}$ and $\widehat{V}^{\pm}$ have simple poles at $s=0$ with residue 1, we pass a pole of order $k$ and obtain
 \begin{displaymath}
\begin{split}
\tilde{M}^{\pm}_j  & =   \frac{C}{(k-1)!}  \log \left( \frac{X^{k-1} \Delta^{-1}}{  (X\Delta)^{(j-1)}}\right)^{k-1}  + O\left((\log X)^{k-2}\right)\\
& =  \frac{C(k-j)^{k-1} }{(k-1)!}  (\log X)^{k-1}  + O\left( (\log X)^{k-2} \log \Delta\right). 
  \end{split}
  \end{displaymath} 
We substitute this back into \eqref{error-F}, combine it with \eqref{inequality}, \eqref{mainterm} and \eqref{finalE} and use   the  formula (cf.\ e.g.\ \cite[0.154.3]{GR}) 
 $$\frac{1}{(k-1)!}\sum_{j=1}^{k-1} (-1)^{j-1}  \left(\begin{array}{l} k\\ j \end{array} \right) (k-j)^{k-1}  = \frac{k^{k-1}}{(k-1)!}$$
to conclude that
\begin{displaymath}
\begin{split}
\mathcal{M}(\mathcal{R}_X) = & \frac{C \cdot \text{vol}(\mathcal{R})}{(k-1)!} X^{k-1} (\log X^k)^{k-1} \\
&+  O_{\varepsilon}\left(X^{k-1} \left( (\log X)^{k-2 + \frac{1}{k-1} +\varepsilon} + (\log X)^{k-2} \log \Delta + (\log X)^{O(1)} \Delta^{-1}\right) \right).
\end{split}
\end{displaymath}
Choosing $\Delta = (\log X)^B$ for a  sufficiently large constant $B$ completes the proof of the theorem.

\end{document}